\documentclass{article}
\usepackage{graphicx} 

\title{Extremal $K_t$-intersecting graph families for large $t$}

\usepackage{amsmath, amssymb, amsthm, url, tikz, verbatim, hyperref}

\newtheorem{theorem}{Theorem}[section]

\newtheorem{lemma}[theorem]{Lemma}

\newtheorem{question}[theorem]{Question}

\newtheorem{conjecture}[theorem]{Conjecture}

\newcommand{\ft}{\lfloor\sqrt{t}\rfloor}
\makeatletter
\NewCommandCopy\@@pmod\pmod
\DeclareRobustCommand{\pmod}{\@ifstar\@pmods\@@pmod}
\def\@pmods#1{\mkern4mu({\operator@font mod}\mkern 6mu#1)}
\makeatother

\usetikzlibrary{calc}
\usetikzlibrary{patterns}
\usetikzlibrary{decorations.markings}
\usetikzlibrary{arrows,shapes.geometric,%
	decorations.pathreplacing,shapes,shadows}

\date{}

\author{}

\author{Yijia Fang\thanks{Department of Mathematics, National University of Singapore. Email: fangyijia@u.nus.edu.}
\and
Hao Huang\thanks{Department of Mathematics, National University of Singapore. Email: huanghao@nus.edu.sg. Research supported in part by a start-up grant at NUS and an MOE Academic Research Fund (AcRF) Tier 1 grant A-8003627.}}
\begin{document}

\maketitle

\begin{abstract}
In this paper, we show that for every sufficiently large positive integer $t$ and every integer $n \ge t$, every $K_t$-intersecting family $\mathcal{G}$ of graphs on vertex set $[n]$ has size at most $2^{\binom{n}{2}-\binom{t}{2}}$. Furthermore, when the equality is attained, $\mathcal{G}$ must consist of all graphs on $[n]$ that contain some fixed copy of $K_t$. Our result settles a conjecture of Ellis, Filmus, and Friedgut for sufficiently large $t$. We also present a short proof for the odd-cycle-intersecting case, slightly different from the existing proofs in the literature.
\end{abstract}
\section{Introduction}
A family $\mathcal{F}$ of subsets of $[n]$ is {\it $r$-intersecting} if any two members of it intersect in at least $r$ elements. A folklore result says every $1$-intersecting family of subsets of $[n]$ has size at most $2^{n-1}$, attained by taking all sets containing a fixed element. This result can be viewed as a non-uniform analogue of the celebrated Erd\H os--Ko--Rado Theorem \cite{EKR61}. One might expect that for $r \ge 2$, the family of all sets containing a fixed $r$-subset of $[n]$, which has size $2^{n-r}$, is also optimal. However, assuming $n+r$ is even, the family of all sets of size at least $(n+r)/2$ is also $r$-intersecting. For fixed $r$ and $n$ tends to infinity, its size is $(1-o_n(1))2^{n-1}$, significantly larger than $2^{n-r}$. Katona \cite{K64} determined the exact maximum size of an $r$-intersecting family, showing in particular that this construction is optimal when $n+r$ is even.

Requiring the pairwise intersection to contain $r$ elements with additional structure could significantly change the nature of this intersection problem. This direction was pioneered by Simonovits and S\'os \cite{SS78, SS80} about half a century ago. A family $\mathcal{G}$ of graphs on vertex set $[n]$ is called {\it triangle-intersecting} if for every two graphs $G_1, G_2 \in \mathcal{G}$, $G_1 \cap G_2$ contains a copy of a triangle. Treating the $\binom{n}{2}$ edges of $K_n$ as elements, such $\mathcal{G}$ is clearly a $3$-intersecting family on a ground set of $\binom{n}{2}$ elements. In 1976, Simonovits and S\'os conjectured that such a family has size at most $2^{\binom{n}{2}-3}$, with equality attained by the family of all graphs containing a fixed triangle. Chung, Frankl, Graham, and Shearer \cite{CFGS86} used the entropy method to prove an upper bound $2^{\binom{n}{2}-2}$. In 2012, Ellis, Filmus, and Friedgut \cite{EFF12} resolved the Simonovits--S\'os Conjecture, and further established multiple strengthenings and extensions using a Fourier-analytic method. Similar to \cite{F08}, this proof can also be viewed as an application of a weighted version of Hoffman's ratio bound \cite{HaemersRatioBound} from spectral graph theory.

It is natural to ask what happens if we replace the triangle by an arbitrary graph $H$. Given a fixed graph $H$, we say that a family $\mathcal{G}$ of graphs on $[n]$ is {\it $H$-intersecting}, if for every two graphs $G_1, G_2 \in \mathcal{G}$, $G_1 \cap G_2$ contains a copy of $H$. The Ellis--Filmus--Friedgut Theorem inspires the following natural conjecture in \cite{EFF12}.

\begin{conjecture}\label{conj_kt}
For positive integers $n \ge t \ge 2$, every $K_t$-intersecting family $\mathcal{G}$ of graphs on $[n]$ has size
$$|\mathcal{G}| \le 2^{\binom{n}{2}-\binom{t}{2}}.$$
\end{conjecture}

This conjecture was confirmed for $t=4$ by Berger and Zhao \cite{BZ23}, and for $t=5$ by Mani and Zhang \cite{MZ25}, both using the Fourier/spectral technique. One might even be tempted to conjecture that regardless of whether $H$ is complete, the maximum size of an $H$-intersecting family $\mathcal{G}$ of graphs on $[n]$ is always $2^{\binom{n}{2}-e(H)}$ as long as $n \ge |V(H)|$, or equivalently, the maximum density $|\mathcal{G}|/2^{\binom{n}{2}}$ of such a family is $2^{-e(H)}$. However, this already fails for stars $H=K_{1,r}$ for $r \ge 2$: one can still apply the Katona construction to edges incident to a fixed vertex, while leaving all other edges unrestricted. For sufficiently large $n$, this gives density $\frac{1}{2}-o(1)$, exceeding $2^{-e(H)}$. Alon observed that for every fixed star-forest $H$ (which includes $K_{1, r}$), the maximum density of an $H$-intersecting family is $\frac{1}{2}-o(1)$. He further conjectured (see \cite{AS}) that there is some universal constant $c>0$ such that whenever $H$ is not a star-forest, every $H$-intersecting family has density at most $\frac{1}{2}-c$. It would suffice to prove this conjecture for the $3$-edge-path $P_4$. Han and Wang \cite{HW26} recently resolved this case, building on techniques developed by Gillott \cite{G26} who showed that every cycle-intersecting family has density bounded away from $\frac{1}{2}$ by an absolute constant.

There have been many other works on the $H$-intersecting problem. Christofides \cite{C12} constructed a $P_4$-intersecting family of density $\frac{17}{128} >\frac{1}{8}=2^{-e(P_4)}$ (see \cite{E22})  beating the canonical construction, disproving a conjecture from \cite{CFGS86}. For complete bipartite graphs, Balogh and Linz \cite{BL24} showed that for $t$ sufficiently large in $s$, there exists a $K_{s,t}$-intersecting family with density strictly greater than $2^{-st}=2^{-e(K_{s,t})}$. Their constructions also extend to complete multipartite graphs with one part sufficiently large compared to others. Sason \cite{S25} used the entropy method to establish an uniform upper bound for all $H$: the density of an $H$-intersecting family is at most $2^{-(\chi(H)-1)}$, extending the results in \cite{CFGS86}.

In this paper, we focus on the case when $H$ is the complete graph $K_t$. We prove Conjecture \ref{conj_kt} for sufficiently large $t$, and also establish the uniqueness of the extremal construction.

\begin{theorem}\label{thm_main}
There exists a positive integer $t_0$, such that for all $n \ge t \ge t_0$, if a family $\mathcal{G}$ of graphs on vertex set $[n]$ is $K_t$-intersecting, then
$$|\mathcal{G}| \le 2^{\binom{n}{2}-\binom{t}{2}}.$$
Furthermore, the equality is attained if and only if $\mathcal{G}$ consists of graphs on $[n]$ containing a fixed copy of $K_t$.
\end{theorem}

The rest of this paper is organized as follows. In the next section, we use the standard Fourier analytic technique developed in \cite{EFF12} to reduce Theorem \ref{thm_main} to Theorem \ref{thm_key}, that is, to finding a function $f$ whose values are restricted and  Fourier coefficients vanish on graphs containing $K_t$. In Section \ref{sec_small_n}, we construct such a function $f$ for $n$ at most polynomial in $t$. Section \ref{sec_interpolation} contains a technical interpolation lemma, which will be used to extend the previous construction of $f$ to all $n$ in Section \ref{sec_large_n}. Section \ref{sec_k3} is somewhat independent from the rest of the paper. There we present a new proof for the odd-cycle-intersecting case (that would imply the $K_3$ result), based mainly on the monotonicity of extremal families. The argument takes a slightly different approach from the existing proofs in the literature and is arguably simpler. In the final section, we discuss the limitation of our method in this paper, potential improvements and generalizations, and some other related  problems to be explored in the future.

\section{The spectral/Fourier reduction}\label{sec_reduction}
Given a vertex set $V$, let $K_V$ be the complete graph on the vertex set $V$. For a given function $f: 2^{E(K_V)} \rightarrow \mathbb{R}$ defined on all spanning subgraphs $G$ of $K_V$, we use the following standard form of Fourier transform:
$$f(S)=\sum_{T \subseteq K_V} \widehat{f}(T) (-1)^{|S \cap T|}, \qquad \widehat{f}(T)=\frac{1}{2^{\binom{|V|}{2}}}\sum_{S \subseteq K_V} f(S) (-1)^{|S \cap T|}.$$
We denote by $\mathcal{F}_q(V)$ the set of all functions $f$ defined on $T \subseteq K_V$, such that  $\widehat{f}$ is supported on graphs with clique number at most $q$ (i.e. $K_{q+1}$-free graphs). When a function $f$ is defined on graphs of every order, then we say $f \in \mathcal{F}_q$ if the restriction of $f$ on every subset $V$ lies in $\mathcal{F}_q(V)$. 

Given a graph $S$, we let $V^+(S) \subseteq V(S)$ be its subset of non-isolated vertices. A function $f$ defined on graphs of every order is called {\it isomorphism-invariant}, if for two isomorphic graphs $S_1 \cong S_2$, we have $f(S_1)=f(S_2)$. $f$ is called {\it isolate-invariant} if whenever appending or removing isolated vertices, the value of $f$ does not change.

\begin{theorem}\label{thm_key}
There exists $t_0$ such that for every $t \ge t_0$, there exists $f \in \mathcal{F}_{t-1}$ that is isomorphism-invariant and isolate-invariant, with\\
(1) $f(\emptyset)=2^{\binom{t}{2}}-1$,\\
(2) $|f(S)| \le 1$ when $|V^+(S)|>0$,\\
(3) $|f(S)| \le t^{-7}$ when $|V^+(S)|>t$.
\end{theorem}

The proof of Theorem \ref{thm_key} will appear in the upcoming sections. In this section, we only show that it would imply our main Theorem \ref{thm_main}. As a matter of fact, if we only want to establish the upper bound $|\mathcal{G}| \le 2^{\binom{n}{2}-\binom{t}{2}}$, it suffices to replace (2) and (3) by a common upper bound $|f(S)| \le 1$ when $S \neq \emptyset$. To further establish the uniqueness of the extremal example, it suffices to replace $t^{-7}$ by any constant strictly smaller than $1$ in (3). Even though the statement of Theorem \ref{thm_key} and its implication of Theorem \ref{thm_main} below are written using the language of Fourier analysis, for readers familiar with spectral graph theory, this is no more than proving Theorem \ref{thm_main} by bounding the independence number of the associated $2^{\binom{n}{2}}$-vertex graph using a weighted version of Hoffman's ratio bound.

\begin{proof}[Proof of Theorem \ref{thm_main}]
We take the $f \in \mathcal{F}_{t-1}$ from Theorem \ref{thm_key} and consider its restriction on the vertex set $[n]$. With slight abuse of notation we still call it $f$. By the definition of $\mathcal{F}_{t-1}$, we know that $\widehat{f}$ vanishes on graphs whose clique number is at least $t$. For two graphs $G, H \in \mathcal{G}$, we have $G \cap H \subseteq \overline{G \Delta H}$, so $\omega(\overline{G \Delta H}) \ge \omega(G \cap H) \ge t$. Therefore, if we define a matrix $A$ whose rows and columns correspond to all graphs on $[n]$, such that $A_{G, H}=\widehat{f}(\overline{G \Delta H})$, then by setting $h=1_{\mathcal{G}}$, clearly we have
$\langle h, Ah \rangle =0$.

We define $\phi_S(G)=(-1)^{|S \cap G|}$ for $S,G \subseteq K_n$. It is not hard to see that 
\begin{align*}
(A\phi_S)(G)&=\sum_{H \subseteq K_n} A_{G, H} \phi_S(H)=\sum_{H \subseteq K_n} \widehat{f}(\overline{G \Delta H}) \phi_S(H)\\
&=\phi_S(G) \sum_{D \subseteq K_n} \widehat{f}(\overline{D}) \phi_S(D)
= \phi_S(G)\cdot (-1)^{|S|}\sum_{D \subseteq K_n}  \widehat{f}(D) \phi_S(D)\\
&= (-1)^{|S|}f(S) \phi_S(G).
\end{align*}
Therefore $\phi_S$ is an eigenvector corresponding to eigenvalue $(-1)^{|S|}f(S)$. It is also straightforward to check that $\{\phi_S\}$ form an orthonormal basis. Therefore, we can write 
$h=\sum_S c_S \phi_S$, and 
\begin{align*}
0 = \langle h, Ah\rangle = \sum_S (-1)^{|S|} f(S) c_S^2.
\end{align*}
Since all $f(S)$ except $S=\emptyset$ satisfy $|f(S)| \le 1$, we have 
$$0 \ge f({\emptyset}) c_{\emptyset}^2 - \sum_{S\neq \emptyset} c_S^2=2^{\binom{t}{2}} c_{\emptyset}^2 - \sum_S c_S^2.$$
We also have 
$$c_{\emptyset}=\langle h, \phi_{\emptyset}\rangle = \frac{\sum_S h(S)}{2^{\binom{n}{2}}}=\frac{|\mathcal{G}|}{2^{\binom{n}{2}}},\qquad \sum_S c_S^2 = \langle h, h\rangle = \frac{|\mathcal{G}|}{2^{\binom{n}{2}}}.$$
Plugged into the previous inequality, we immediately have
$$|\mathcal{G}| \le 2^{\binom{n}{2}-\binom{t}{2}}.$$

For the remainder of this proof we investigate when the equality can be attained. When that happens, we have that for every $S \neq \emptyset$, $c_S=0$ whenever $(-1)^{|S|}f(S)>-1$, which holds whenever $|V^+(S)|>t$. In other words, $\widehat{h}(S)=0$ when $|V^+(S)|>t$.

We take an inclusion-minimal graph $M \in \mathcal{G}$. Thus $h(M)=1$, and  $h(S)=0$ for every proper subgraph $S \subsetneq M$. We can view $h$ as a function in the variables $x_e$ with $e \in E(K_n)$, then its Fourier expansion is
$$h(x)=\sum_S \widehat{h}(S) \prod_{e \in S} (1-2x_e).$$
Note that whenever we set all coordinates of $x$ outside those from $M$ to be zero, the value of $h$ is equal to the product of $x_e$ for all $e \in M$. Therefore we have
$$\prod_{e \in M} x_e = \sum_S \widehat{h}(S) \prod_{e \in S \cap M} (1-2x_e).$$
Comparing the coefficients of $\prod_{e \in M} x_e$ from the two sides, we have
$$1=(-2)^{|M|} \sum_{S: M \subseteq S} \widehat{h}(S).$$
In particular there is some $S \supseteq M$ with $\widehat{h}(S) \neq 0$.  Since $\widehat{h}$ vanishes when the number of non-isolated vertices is strictly greater than $t$, we have that $S$ has at most $t$ non-isolated vertices. Thus $M$ also has at most $t$ non-isolated vertices. Since $\mathcal{G}$ is $K_t$-intersecting, in particular $\omega(M) \ge t$. The only possibility is $M \cong K_t$, ignoring isolated vertices. Since for every other graph $G \in \mathcal{G}$, $\omega(G \cap M) \ge t$, this immediately gives $G \supseteq M$. The number of graphs containing $M$ is precisely $2^{\binom{n}{2}-\binom{t}{2}}$, therefore $\mathcal{G}$ consists of all such graphs. This completes the proof of the uniqueness part of Theorem \ref{thm_main}.
\end{proof}

\section{Constructing the function $f$ for $t \le n \le t^{16}$}\label{sec_small_n}

The goal of this section is to prove the existence of such a function $f$ for $t \le n \le t^{16}$. In this interval, for (3) we will get a slightly stronger upper bound $|f(S)| \le t^{-10}$. Such $f$ will be the cornerstone for us to extend our constructions of $f$ to larger $n$ in later sections.

\begin{theorem}\label{thm_small_n}
There exists $t_0$ such that for every positive integer $t \ge t_0$, and every integer $t \le n \le t^{16}$, there exists an isomorphism-invariant function $f \in \mathcal{F}_{t-1}([n])$ such that \\
(1) $f(\emptyset)=2^{\binom{t}{2}}-1$,\\
(2) $f(S) = (-1)^{e(S)+1}$ for $S$ with $0<|V^+(S)| \le t$,\\
(3) $|f(S)| \le t^{-10}$ for $S$ with $|V^+(S)|>t$.
\end{theorem}

Before turning to the details, we outline the main ideas underlying the proof. We will construct a function $g: 2^{E(K_n)} \rightarrow \mathbb{R}$, such that $f$ defined by $f(S)=(-1)^{e(S)}g(S)$ satisfies the requirements from Theorem \ref{thm_small_n}. To guarantee that $f \in \mathcal{F}_{t-1}([n])$, we introduce the following projection operator $P$: given a function $h: 2^{E(K_n)} \rightarrow \mathbb{R}$ with Fourier expansion
$$h(S)=\sum_{T \subseteq K_n} \widehat{h}(T) (-1)^{|S\cap T|},$$ we define a projection operator $P$ that maps $h$ to a function $Ph: 2^{E(K_n)} \rightarrow \mathbb{R}$ such that 
\begin{equation}\label{eq_p}
(Ph)(S)=\sum_{T \subseteq K_n: \omega(\overline{T}) \ge t} \widehat{h}(T) (-1)^{|S \cap T|}.
\end{equation}
In other words, we only keep those frequencies of graphs whose complement contains a copy of $K_t$, or, $\widehat{Ph}$ is supported on those $T$ with $\omega(\overline{T}) \ge t$. It is also not hard to see that for an arbitrary function $h: 2^{E(K_n)} \rightarrow \mathbb{R}$, if we set $$f(S)=(-1)^{e(S)}(h(S)-(Ph)(S)),$$ then $\widehat{f}$ is supported on $K_t$-free graphs, and thus such $f$ satisfies $f \in \mathcal{F}_{t-1}([n])$. We also obviously have $P^2 = P$. From these observations, it is straightforward to see that
$$f(S)=(-1)^{e(S)} g(S) \in \mathcal{F}_{t-1}([n])~~~\Longleftrightarrow~~~Pg=0.$$

Next we explain how to construct a satisfiable $g$. We start with the base case $n=t$. For this case, we can simply take $g(\emptyset)=2^{\binom{t}{2}}-1$, $g(S)=-1$ when $S \neq \emptyset$. This translates to $f(\emptyset)=2^{\binom{t}{2}}-1$, and $f(S)=(-1)^{e(S)+1}$ when $S \neq \emptyset$. It is not hard to see $\widehat{f}(S)=1$ if $S \neq K_t$, and $0$ when $S=K_t$. So $f \in \mathcal{F}_{t-1}([n])$.

For general $n \ge t$, we start with essentially the same function $g_0$ such that 
\begin{equation}\label{eq_g0}
g_0(\emptyset)=2^{\binom{t}{2}}-1, \qquad g_0(S)=-1~\textup{for~}0<|V^+(S)| \le t, \qquad g_0(S)=0\textup{~for~}|V^+(S)|>t.
\end{equation}
We know $g_0-Pg_0$ would give a new function $g_1$ with $Pg_1=0$. However, when modifying some Fourier coefficients for $g_0$, this might change the value of this function for graphs $S$ with $|V^+(S)| \le t$, for which we need to meet requirements (1) and (2) in Theorem \ref{thm_small_n}. To settle this issue, we can instead subtract only the correction outside of those graphs. To make our argument more precise, throughout this section we let 
$$\mathcal{B}=\{S \subseteq K_n: |V^+(S)| \le t\}, \qquad \mathcal{D}=\{S \subseteq K_n: |V^+(S)| >t\},$$
and define $1_{\mathcal{B}}$ and $1_{\mathcal{D}}$ to be the operators keeping the value of a function at those inputs and changing the rest to be zero. It is not hard to see that $1_{\mathcal{B}}+1_{\mathcal{D}}=I$, the identity operator. So instead of subtracting $Pg_0$ from $g_0$, a natural idea is to subtract $1_{\mathcal{D}}y$ from $g_0$ for some function $y$, since this would not change the value of $g_0$ on $\mathcal{B}$. Recall that we also need $P(g_0-1_{\mathcal{D}}y)=0$. Set $T=P 1_{\mathcal{B}}$. Using $1_{\mathcal{D}}=I-1_{\mathcal{B}}$, this becomes
$$Pg_0 = P(I-1_{\mathcal{B}}) y=(P -T)y.$$
Since $P$ is an orthogonal projection with $P^2=P$, intuitively $P$ is an identity map on the projected subspace. The natural choice for $y$ is thus $(P-T)^{-1} Pg_0=(I-T)^{-1} Pg_0$. At this moment it is not clear why the inverse of $(I-T)$ is well-defined. Later in this section, we will prove that $T$ is a strict contraction in an invariant subspace $\mathcal{C}$ containing $Pg_0$, and thus now we can instead define
$$y=\sum_{i \ge 0} T^i Pg_0,$$
and let
$$g=g_0-1_{\mathcal{D}} y=g_0-1_{\mathcal{D}} \sum_{i \ge 0} T^i Pg_0.$$
Let us carefully verify $Pg=0$ in the proof below. Note that in this proof we have not shown that $\sum_{i \ge 0} T^i Pg_0$ converges. This will become self-evident later when we prove the operator norm of $T$ on that subspace $\mathcal{C}$ is strictly less than $1$ in Theorem \ref{thm_final_estimate}.
\begin{lemma}\label{lem_f_in_F}
For the function $g$ constructed above, let $f(S)=(-1)^{e(S)} g(S)$, then we have
$$f \in \mathcal{F}_{t-1}([n]).$$
\end{lemma}
\begin{proof}
First, since $PT=P^2 1_{\mathcal{B}}=P1_{\mathcal{B}}=T$, we have
\begin{align*}
Pg &=  Pg_0 - P 1_{\mathcal{D}} y = Pg_0 - P (I-1_{\mathcal{B}}) \sum_{i \ge 0} T^i P g_0\\
&= Pg_0 -  \sum_{i \ge 0} \left(P T^i P g_0 -P1_{\mathcal{B}} T^i P g_0\right)= Pg_0 -  \sum_{i \ge 0} \left(T^i P g_0 -T^{i+1} P g_0\right)\\
&=Pg_0 -P g_0 = 0.
\end{align*}
Therefore, $f(S)=(-1)^{e(S)}(g(S)-(Pg)(S)) = (-1)^{e(S)} g(S) \in \mathcal{F}_{t-1}([n])$. 
\end{proof}

By the definition of $g$, clearly for $S \in \mathcal{B}$, $g(S)=g_0(S)$. These would give (1) and (2) for $f(S)=(-1)^{e(S)}g(S)$. So it remains to prove $|f(S)|=|g(S)| \le t^{-10}$ for $S \in \mathcal{D}$. Note that for $S \in \mathcal{D}$, $g(S)=g_0(S)-1_{\mathcal{D}}y(S)=-y(S)$. So it suffices to bound $|y(S)|$ by $t^{-10}$.

To prove this, we introduce the following weighted norm for an arbitrary function $u$. We choose the weight $w(S)$ for a graph $S$ as follows:
$$w(S)=\begin{cases}
2^{\binom{t}{2}}, \qquad & S=\emptyset,\\
n^{t-|V^+(S)|},\qquad & S \in \mathcal{B} \setminus \{\emptyset\},\\
1,\qquad & S \in \mathcal{D}.
\end{cases}$$
and define the {\it $w$-norm} of a function $u$ as
$$\|u\|_w=\max_S \frac{|u(S)|}{w(S)}.$$
We will prove the following three key results:

\begin{itemize}
\item Define $$\mathcal{C}=\left\{u: \sum_{S \subseteq K_V} u(S)=0 \textup{~for~every~} V \subseteq [n], |V|=t\right\}.$$
Then for any $u \in \mathcal{C}$, $Tu \in \mathcal{C}$.
\item For $g_0$ we defined earlier in \eqref{eq_g0}, we have $Pg_0 \in \mathcal{C}$, and
\begin{equation}\label{ineq_norm_pg0}
\|Pg_0\|_w \le \frac{1}{2}t^{-10}.
\end{equation}
\item 
For every $u \in \mathcal{C}$, we have

\begin{equation}\label{ineq_norm_t}
\|Tu\|_w \le \frac{1}{2}t^{-10} \|u\|_w.
\end{equation}
\end{itemize}

We postpone the proof of these claims to Lemma \ref{lem_C} and Theorem \ref{thm_final_estimate} later on. Assuming all these claims, now we can complete the proof of Theorem \ref{thm_small_n}.

\begin{proof}[Proof of Theorem \ref{thm_small_n}]
From \eqref{ineq_norm_pg0}, we have $\|Pg_0\|_w \le \frac{1}{2}t^{-10}$. Then by repeatedly applying \eqref{ineq_norm_t} and using the fact that $Pg_0 \in \mathcal{C}$ and $T$ maps $\mathcal{C}$ to itself, we can use induction on $i$ to prove
$$\|T^i Pg_0\|_w \le \frac{1}{2} t^{-10} \|T^{i-1}Pg_0\|_w \le \cdots \le \left(\frac{1}{2}t^{-10}\right)^{i+1}.$$
Therefore by the triangle inequality, for large $t$,
$$\|y\|_w=\left\|\sum_{i \ge 0} T^i P g_0 \right\|_w \le \sum_{i \ge 0} \left(\frac{1}{2}t^{-10}\right)^{i+1} = \frac{\frac{1}{2}t^{-10}}{1-\frac{1}{2}t^{-10}}<t^{-10}.$$
In other words, when $t$ is sufficiently large, for $S \in \mathcal{D}$ we have
$$y(S) \le \|y\|_w w(S) =\|y\|_w <t^{-10}.$$
This implies for $S \in \mathcal{D}$, 
$$|f(S)|=|g(S)|=|y(S)|<t^{-10},$$
So (3) is verified. For $S=\emptyset$ and $S \in \mathcal{B}\setminus \{\emptyset\}$, clearly we have $g(S)=g_0(S)$ which cause $f(S)$ to satisfy the requirements (1) and (2). Lemma \ref{lem_f_in_F} shows that $f \in \mathcal{F}_{t-1}([n])$. These facts altogether complete the proof.
\end{proof}

The remainder of this section is devoted to proving the earlier claims.

\begin{lemma}\label{lem_C}
We have $Pg_0 \in \mathcal{C}$. Also, for $u \in \mathcal{C}$, we have $T u \in  \mathcal{C}$.
\end{lemma}
\begin{proof}
By our choice of $g_0$, clearly $g_0 \in \mathcal{C}$. For a $t$-set $V$ of $[n]$, we define $\phi_V(S)=1_{S \subseteq K_V}$. Then $u \in \mathcal{C}$ if and only if $\langle u, \phi_V\rangle=0$ for every $|V|=t$. The Fourier transform of $\phi_V$ is supported on graphs $H$ with no edges inside $V$. Since $P$ only changes those Fourier coefficients for graphs $H$ with $\alpha(H)=\omega(\overline{H}) <t$ and $|V|=t$, so we have $P\phi_V=\phi_V$. Since $P$ is an orthogonal projection, we have 
$$\langle Pg_0, \phi_V \rangle = \langle g_0, P\phi_V\rangle =\langle g_0, \phi_V\rangle=0,$$
thus $Pg_0 \in \mathcal{C}$.

Next we show for $u \in \mathcal{C}$, $Tu$ is still in $\mathcal{C}$. The key observation here is: $1_{\mathcal{B}}\phi_V=\phi_V$, since $|V|=t$ and thus any non-vanishing $\phi_V(S)$ has at most $t$ non-isolated vertices. Therefore, for every $|V|=t$,
$$\langle Tu, \phi_V\rangle =\langle P1_\mathcal{B} u, \phi_V\rangle = \langle 1_\mathcal{B} u, P \phi_V\rangle=\langle 1_\mathcal{B} u,\phi_V\rangle=\langle  u, 1_\mathcal{B} \phi_V\rangle =\langle u, \phi_V\rangle=0.$$
So we also have $Tu \in \mathcal{C}$.
\end{proof}

For each subset $W \subseteq [n]$, we define a function $a_W$ on graphs $G$ on vertex set $W$, with
$$a_W(G)=\sum_{U: U \subseteq W} (-1)^{|W|-|U|} \cdot 1_{\{\omega(G[U]) \ge t\}}.$$
The motivation for defining such $a_W$ is such that they satisfy for every graph $G$ on $V \subseteq [n]$,
\begin{equation}\label{eq_aw_sum}
1_{\{\omega(G) \ge t\}}=\sum_{W \subseteq V} a_W(G[W]).
\end{equation}
We let $A_W$ act on the Fourier coefficients by 
\begin{equation}\label{eq_Aw}
\widehat{A_W u}(H)=a_W(\overline{H[W]}) \widehat{u}(H).
\end{equation}
Then the following lemma decomposes the projection operator $P$ into operators $A_W$'s labeled by subsets $W$ of $[n]$.
\begin{lemma}\label{lem_p_decompose}
We have
$$P=\sum_{W \subseteq [n]} A_W.$$
Furthermore, for every function $u$ and every graph $S$, 
$$(A_W u) (S)=\sum_{T \subseteq K_n} u(T) \cdot \mathbb{E}_H\left[a_W(\overline{H[W]})(-1)^{|(S \Delta T) \cap H|}\right],$$
where $H$ is a uniform random subgraph on $[n]$.
\end{lemma}
\begin{proof}
Recall that our projection operator $P$ first defined in \eqref{eq_p} was chosen such that
$$\widehat{Pu}(H)=\widehat{u}(H) \cdot 1_{\{\omega(\overline{H}) \ge t\}}.$$
Using \eqref{eq_Aw} and \eqref{eq_aw_sum}, we have that for a subgraph $H$ on vertex set $[n]$,
\begin{align*}
\widehat{Pu}(H)&=\widehat{u}(H) \cdot 1_{\{\omega(\overline{H}) \ge t\}}
=\widehat{u}(H) \sum_{W \subseteq [n]} a_W(\overline{H[W]})
=\sum_{W \subseteq [n]} \widehat{A_W u}(H).
\end{align*}
This gives the decomposition $P=\sum_{W \subseteq [n]} A_W$. 

To prove the second assertion, we plug the Fourier inversion formula for $\widehat{u}(H)$ into \eqref{eq_Aw}, we get
\begin{align*}
\widehat{A_W u}(H)&=a_W(\overline{H[W]}) \cdot \frac{1}{2^{\binom{n}{2}}} \sum_{T \subseteq K_n} u(T) (-1)^{|T \cap H|}.
\end{align*}
Therefore
\begin{align*}
(A_W u)(S)&=\sum_{H \subseteq K_n} \widehat{A_W u}(H) (-1)^{|H \cap S|}\\
&=\sum_{T \subseteq K_n} u(T) \left[ 2^{-\binom{n}{2}}\sum_{H \subseteq K_n}\ a_W(\overline{H[W]})(-1)^{|T \cap H|}(-1)^{|S \cap H|}\right]\\
&=\sum_{T \subseteq K_n} u(T) \left[ 2^{-\binom{n}{2}}\sum_{H \subseteq K_n}\ a_W(\overline{H[W]})(-1)^{|(S \Delta T) \cap H|}\right]\\
&=\sum_{T \subseteq K_n} u(T) \cdot \mathbb{E}_H\left[a_W(\overline{H[W]})(-1)^{|(S \Delta T) \cap H|}\right],
\end{align*}
establishing the second assertion.
\end{proof}

Our next lemma is a technical estimate on the weights $w(S)$ that we have chosen earlier. We leave its proof to the appendix.
\begin{lemma} \label{lem_pair}
For sufficient large $t$, $t \le n \le t^{16}$, and every integer $r \ge 0$ and every graph $S$,
$$\frac{1}{w(S)} \sum_{W \subseteq [n]: |W|=t+r} \sum_{T \in \mathcal{B}: S \Delta T \subseteq K_W} w(T) \le 4n^{2r}2^{\binom{t}{2}}.$$
Moreover, if $r=0$ and the outer sum is restricted to sets $W$ for which $S \setminus K_W \neq \emptyset$, then the right hand side improves to $4n^3 \cdot 2^{(t^2-2t+1)/2}$.
\end{lemma}

Assuming Lemma \ref{lem_pair}, the next lemma uses $a_W$ to bound the sum of $A_W$ at a given level.

\begin{lemma}\label{lem_bound_AW_by_aW}
For every $r \ge 0$ satisfying $r+t \le n$, we have that for every function $u$,
$$\left\|\sum_{|W|=t+r} A_W(1_{\mathcal{B}} u)\right\|_w \le 4 \cdot 2^{\binom{t}{2}} n^{2r} \|u\|_w \cdot \max_{|W|=t+r} \mathbb{E}_G|a_W(G)|.$$
Here $\mathbb{E}_G|a_W(G)|$ is the expected value of $|a_W(G)|$ for a uniform random graph $G$ on $W$.
\end{lemma}
\begin{proof}
By Lemma \ref{lem_p_decompose}, for fixed $W$ and each $S$ we have
\begin{align*}
(A_W (1_{\mathcal{B}}u))(S) &= \sum_{T \in \mathcal{B}} u(T) \cdot  \mathbb{E}_H\left[a_W(\overline{H[W]})(-1)^{|(S \Delta T)\cap H|}\right].
\end{align*}
The expectation inside the sum on the right hand side vanishes unless $S \Delta T \subseteq K_W$ since adding or deleting an edge outside $K_W$ would not change $a_W(\overline{H[W]})$ but would change the sign of $(-1)^{|(S \Delta T) \cap H|}$. Therefore 
\begin{align*}
(A_W (1_{\mathcal{B}}u))(S) &= \sum_{T \in \mathcal{B}: S \Delta T \subseteq K_W} u(T) \cdot  \mathbb{E}_H\left[a_W(\overline{H[W]})(-1)^{|(S \Delta T)\cap H|}\right].\\
& \le  \mathbb{E}_G|a_W(G)| \cdot \left(\sum_{T \in \mathcal{B}: S \Delta T \subseteq K_W} |u(T)| \right)\\
&\le \mathbb{E}_G|a_W(G)| \cdot \left(\sum_{T \in \mathcal{B}: S \Delta T \subseteq K_W} w(T) \right) \|u\|_w.
\end{align*}
Adding this inequality for all subsets $W$ of size $t+r$ and divide by $w(S)$ we have
\begin{align*}
\frac{|\sum_{|W|=t+r} A_W (1_{\mathcal{B}}u))(S)|}{w(S)} & \le  \max_{|W|=t+r}\mathbb{E}_G|a_W(G)|\cdot \frac{1}{w(S)} \cdot \left(\sum_{|W|=t+r} \sum_{T \in \mathcal{B}: S \Delta T \subseteq K_W} w(T) \right) \|u\|_w.
\end{align*}
Applying Lemma \ref{lem_pair} immediately concludes the proof.
\end{proof}

Note that we have 
$$P=\sum_W A_W= \left(\sum_{|W|\le t} A_W\right)  + \left(\sum_{t<|W|\le t+\ft} A_W\right) + \left(\sum_{|W|>t+\ft} A_W\right).$$
In the next three lemmas we will estimate the $w$-norm of the operator $T=P1_{\mathcal{B}}$ one by one, based on this decomposition.

\begin{lemma}[Estimate for $|W| \le t$]\label{lem_small}~\\
For $|W|<t$, $A_W=0$. And for $|W|=t$, 
$$(A_W u)(S)=2^{-\binom{t}{2}}\sum_{T \subseteq K_W} u(S \Delta T).$$
As a consequence, for function $u \in \mathcal{C}$, we have
$$\left\|\sum_{|W|=t} A_W (1_{\mathcal{B}} u)\right\|_w \le \gamma_1 \|u\|_w$$
for $$\gamma_1=4n^3 2^{-(t-1)/2}.$$
\end{lemma}
\begin{proof}
For $|W|<t$, it is easy to see from the definition of $a_W$ that in the sum, $G[U]$ has at most $|W|<t$ vertices, so it is not possible for $G[U]$ to contain a $t$-clique. Therefore $a_W=0$, and thus $A_W=0$ from \eqref{eq_Aw}.

Now assume $|W|=t$. In this case $a_W(\overline{H[W]})=1_{\{\overline{H[W]}\cong K_t\}}=1_{\{H \cap K_W=\emptyset\}}$ By \eqref{eq_Aw}, we have
$\widehat{A_W u}(H)=\widehat{u}(H) \cdot 1_{\{H \cap K_W=\emptyset\}}$. Thus using the Fourier expansion of $A_W u$,
$$(A_W u)(S)=\sum_{H \subseteq K_n} \widehat{u}(H) \cdot 1_{\{H \cap K_W=\emptyset\}}\cdot (-1)^{|S \cap H|}.$$
We observe that the following equality holds:
$$1_{\{H \cap K_W=\emptyset\}}=2^{-\binom{t}{2}} \sum_{T \subseteq K_W} (-1)^{|T \cap H|},$$
since for every $H$ disjoint from $K_W$, there are $2^{\binom{t}{2}}$ subsets $T \subseteq K_W$ that give $1$ in the sum, and once $H \cap K_W \neq \emptyset$ the right hand side becomes $0$ since the $1$'s and $-1$'s cancel out. Using this expression, we have
\begin{align*}
(A_W u)(S)&=2^{-\binom{t}{2}}\sum_{T \subseteq K_W}\sum_{H \subseteq K_n} \widehat{u}(H) \cdot (-1)^{|S \cap H|} (-1)^{|T \cap H|}\\
&=2^{-\binom{t}{2}}\sum_{T \subseteq K_W}\sum_{H \subseteq K_n} \widehat{u}(H) \cdot (-1)^{|(S \Delta T) \cap H|} =2^{-\binom{t}{2}}\sum_{T \subseteq K_W}u(S \Delta T),
\end{align*}
completing the proof of the identity.

For the estimates, first observe that 
$$A_W(1_{\mathcal{B}} u)(S)=2^{-\binom{t}{2}} \sum_{T \subseteq K_W} 1_{\mathcal{B}} u(S \Delta T).$$
If $S \subseteq K_W$, then $S \Delta T \subseteq K_W$ and thus $1_{\mathcal{B}} u (S \Delta T)=u(S \Delta T)$ which ranges over all subgraphs $S \Delta T$ of $K_W$. Therefore in this case $$A_W(1_{\mathcal{B}} u)(S)=2^{-\binom{t}{2}} \sum_{T \subseteq K_W} u(T)=0$$
since $u \in \mathcal{C}$.

Therefore in $\sum_{|W|=t} A_W (1_{\mathcal{B}}u)$, only sets $W \subseteq [n]$ with $S \setminus K_W \neq \emptyset$ can contribute. For such $W$, we have
\begin{align*}
\left|\sum_{|W|=t} A_W (1_{\mathcal{B}}u)(S)\right| &\le  2^{-\binom{t}{2}} \left|\sum_{|W|=t, S\setminus K_W \neq \emptyset} \sum_{V \in \mathcal{B}, V \Delta S \subseteq K_W} u(V)\right|\\
&\le 2^{-\binom{t}{2}} \|u\|_w \sum_{|W|=t, S\setminus K_W \neq \emptyset} \sum_{V \in \mathcal{B}, V \Delta S \subseteq K_W} w(V)\\
&\le  2^{-\binom{t}{2}} \|u\|_w \cdot  4n^3 \cdot 2^{(t^2-2t+1)/2}\cdot w(S)\\
&=\gamma_1 \|u\|_w \cdot w(S). 
\end{align*}
The last inequality follows from the second part of Lemma \ref{lem_pair}. This immediately implies the $w$-norm is bounded by $\gamma_1 \|u\|_w$.

\end{proof}

\begin{lemma}[Estimate for $t<|W| \le t+\ft$]\label{lem_mid}~\\
For $|W|=t+r$ with $1 \le r \le \ft$ and $r+t \le n$, we have that for a uniform random graph $G$ on $W$,
$$\mathbb{E}_G|a_W(G)| \le 2(r+1)n^{3r}2^{-\binom{t}{2}-t}.$$
Consequently, for every function $u$,
$$\left\|\sum_{t<|W|\le t+\ft} A_W (1_{\mathcal{B}} u)\right\|_w \le \gamma_2\|u\|_w$$
for 
$$\gamma_2=10tn^{5\sqrt{t}} 2^{-t}.$$
\end{lemma}

\begin{proof}
We start by estimating $\mathbb{E}_G|a_W(G)|$.
Recall that we defined $a_W(G)=\sum_{U \subseteq W} (-1)^{|W|-|U|} 1_{\{\omega(G[U]) \ge t\}}$. If $G$ has no $t$-clique, clearly $a_W(G)=0$. If $G$ has only one $t$-clique, then $1_{\{\omega(G[U]) \ge t\}}$ is $1$ when the unique clique is contained in $U$ and vanishes otherwise. Since $|W|>t$, we have the parities of $(-1)^{|W|-|U|}$ cancel out, so we still have $a_W(G)=0$ when $G$ has only one $t$-clique. So we only need to consider those $G$ that contain at least two distinct $t$-cliques. These two $t$-cliques share at most $t-1$ vertices, and thus their edge union has at least 
$$2\binom{t}{2}-\binom{t-1}{2}=\binom{t}{2}+(t-1)$$
edges. For a fixed pair of two $t$-cliques, the probability that all these edges are in $G$ is at most $2^{-\binom{t}{2}-(t-1)}$. There are at most $\binom{|W|}{t}^2 = \binom{t+r}{t}^2 = \binom{t+r}{r}^2 \le n^{2r}$ such pairs. The union bound gives 
$$\mathbb{P}[a_W \neq 0] \le n^{2r} 2^{-\binom{t}{2}-(t-1)}.$$
In the expression $a_W(G)=\sum_{U \subseteq W} (-1)^{|W|-|U|} 1_{\{\omega(G[U]) \ge t\}}$, only subsets $U$ of size $|U| \ge t$ can contribute and each of them contributes at most $1$. So 
$$|a_W(G)| \le \sum_{j=t}^{t+r}\binom{t+r}{j}=\sum_{j=0}^r \binom{t+r}{j} \le (r+1)n^r.$$
Multiplying the two bounds together, we have
$$\mathbb{E}_G|a_W(G)| \le 2(r+1)n^{3r}2^{-\binom{t}{2}-t}.$$
Applying Lemma \ref{lem_bound_AW_by_aW}, we immediately have for one fixed $1 \le r \le \sqrt{t}$,
\begin{align*}
\left\|\sum_{|W|=t+r} A_W(1_{\mathcal{B}} u)\right\|_w &\le 4 \cdot 2^{\binom{t}{2}} n^{2r} \cdot \mathbb{E}_G|a_W(G)| \cdot \|u\|_w\le 4 \cdot 2^{\binom{t}{2}} n^{2r} \cdot 2(r+1)n^{3r}2^{-\binom{t}{2}-t} \cdot \|u\|_w\\
&=8(r+1) n^{5r} 2^{-t} \|u\|_w.
\end{align*}
Therefore for large $t$, the triangle inequality gives
\begin{align*}
\left\|\sum_{t<|W|\le t+\ft} A_W(1_{\mathcal{B}} u)\right\|_w &\le \sum_{r=1}^{\ft} 8(r+1)n^{5r}2^{-t}\|u\|_w \le \sqrt{t}\cdot 8(\sqrt{t}+1) n^{5 \ft} 2^{-t} \|u\|_w\\
&\le 10t n^{5\sqrt{t}} 2^{-t} \|u\|_w =\gamma_2 \|u\|_w.
\end{align*}
\end{proof}

\begin{lemma}[Estimate for $|W| > t+\ft$]\label{lem_large}~\\
There exists a constant $c>0$ such that for every sufficiently large $t$, and every $t \le n \le t^{16}$, 
$$\left\|\sum_{|W|> t+\ft} A_W (1_{\mathcal{B}} u)\right\|_w \le \gamma_3\|u\|_w$$
for $$\gamma_3=2^{-ct^{3/2}}.$$
\end{lemma}

\begin{proof}
We denote by $M$ the operator $\sum_{|W|>t+\ft} A_W$. Using $\widehat{A_W u}(H)=a_W(\overline{H[W]}) \widehat{u}(H)$, and $P=\sum_W A_W,$ we have
\begin{equation}\label{eq_M}
\widehat{M u}(H)=\left(\sum_{|W|>t+\ft} a_W(\overline{H[W]})\right) \widehat{u}(H)=: \rho(H) \widehat{u}(H).
\end{equation}
From the previous proof, we learned that $a_W(G)$=0 if $G$ contains at most one $t$-clique. Furthermore, if $G$ has a vertex $v$ not contained in any $t$-clique, in the expression $a_W(G)=\sum_{U \subseteq W} (-1)^{|W|-|U|} 1_{\{\omega(G[U]) \ge t\}}$, we can pair $U$ and $U \cup \{v\}$, then their corresponding terms cancel and this gives  $a_W(G)=0$ as well. So whenever $a_W(G) \neq 0$, every vertex of $W$ belongs to some $t$-clique in $G$.

Suppose for $H$, the Fourier multiplier $\rho(H)$ is non-zero, then for some $|W|>t+\ft$, $a_W(\overline{H[W]}) \neq 0$. Let $G=\overline{H[W]}$. Then $a_W(G) \neq 0$ implies there are at least two $t$-cliques in $G$. Either one of the following two cases must be true:\medskip\\
(1) there exist two $t$-cliques $U, V$ such that $|V \setminus U| \ge \ft$; \smallskip\\
(2) for every pair of $t$-cliques $U \neq V$ in $G$, $|V \setminus U|<\ft$.\medskip\\
Next let us estimate the probability for each case in a uniformly random $H$ (thus $G=\overline{H}$ is also uniformly random).

For the first case, we have $|U \cap V|=t-|V \setminus U| \le t-\ft$. Then $U \cup V$ has at least $2\binom{t}{2}-\binom{t-\ft}{2}$ distinct edges to be present. There are at most $\binom{n}{t}^2$ such pairs $(U, V)$, so a union bound gives
$$\mathbb{P}[\,\textup{the~first~case~occurs}\,]\le 2^{-2\binom{t}{2}+\binom{t-\ft}{2}}\cdot \binom{n}{t}^2\le 2^{-\binom{t}{2}-t\sqrt{t}+O(t\log t)}.$$
The last estimate comes from $n \le t^{16}$.

For the second case, we fix a $t$-clique $U$. Then every $t$-clique $V \neq U$ satisfies $|V \setminus U|<\ft$. For every vertex $v \in W \setminus U$, from the previous observation it must belong to some $t$-clique $V$ in $G$. Thus $v$ is adjacent to every vertex in $V \cap U$, meaning that such $v$ has at most $|U|-|V \cap U|=|V|-|V \cap U|=|V \setminus U|<\ft$ non-neighbors in $U$. And there are at least $|W\setminus U|=|W|-|U|>t+\ft-t=\ft$ choices for such vertex $v$. For a uniform random graph $G$ on $W$, the probability that $U \subseteq G$ is $2^{-\binom{t}{2}}$. The probability that $v \not\in U$ has at most $\ft$ non-neighbors in $U$ is 
$$2^{-t} \sum_{j=0}^{\ft} \binom{t}{j} \le 2^{-t}(\sqrt{t+1}+1)n^{\sqrt{t}} \le n^{\sqrt{t}+1} \cdot 2^{-t}.$$
For different vertices $v \not\in U$, these events are independent. There are at most $\binom{n}{t}$ choices of $U$, and $n^{\sqrt{t}}$ choices of satisfiable vertices $v$. So 
$$\mathbb{P}[\,\textup{the~second~case~occurs}\,]\le \binom{n}{t} n^{\sqrt{t}} \cdot 2^{-\binom{t}{2}} \left(n^{\sqrt{t}+1} 2^{-t}\right)^{\ft} \le 2^{-\binom{t}{2}-t\sqrt{t}+O(t \log t)},$$
again implicitly using $n \le t^{16}$ and the fact that a linear function in $t$ can be absorbed into the $O(t \log t)$ term.

Combining these two cases, we have
$$\mathbb{P}_H(\rho(H) \neq 0) \le 2^{-\binom{t}{2}-t\sqrt{t}+O(t \log t)}.$$
We also have $\rho(H)=1_{\{\omega(\overline{H}) \ge t\}}-\sum_{|W| \le t+\ft} a_W(\overline{H[W]})$, thus we use the same estimate from the proof  in Lemma \ref{lem_mid} and inequality $\binom{n}{t+r}\le \binom{n}{t}n^r$,
$$|\rho(H)| \le 1+\sum_{r=0}^{\ft} \binom{n}{t+r}(r+1)n^r \le 2(\sqrt{t}+1)^2 \binom{n}{t} n^{2\sqrt{t}} \le 2^{O(t \log t)}.$$
Combining these two inequalities, we have that from \eqref{eq_M}, 
\begin{align*}
|M(1_{\mathcal{B}} u)(S)|&=\left|\sum_{H \subseteq K_n} \rho(H) \widehat{1_{\mathcal{B}} u}(H) (-1)^{|S \cap H|}\right|\\
&=2^{-\binom{n}{2}} \left|\sum_{H \subseteq K_n}  \sum_{T \in \mathcal{B}}\rho(H) u(T) (-1)^{|(S\Delta T) \cap H|}\right|\\
&\le \left(\sum_{T \in \mathcal{B}} |u(T)|\right)\cdot \mathbb{P}_H[\rho(H) \neq 0] \cdot \max_H |\rho(H)|. 
\end{align*}
To bound the first sum, we use 
$\sum_{T \in \mathcal{B}} |u(T)| \le \left(\sum_{T \in \mathcal{B}} w(T)\right)\|u\|_w.$ By the definition of $w$, $w(\emptyset)=2^{\binom{t}{2}}$, and $w(T)=n^{t-k}$ for $T \in \mathcal{B}$ with $k$ non-isolated vertices. Note that for a non-empty graph with exactly $n$ vertices and $k$ non-isolated vertices, there are at most $\binom{n}{k}$ choices for these vertices, and $2^{\binom{k}{2}}$ vertices for the adjacency between them. A quick calculation using $\binom{n}{k} \le n^k$ gives
$$\sum_{T \in \mathcal{B}} w(T) \le 2^{\binom{t}{2}} + \sum_{k=1}^t \binom{n}{k}2^{\binom{k}{2}} n^{t-k} \le (t+1)n^t 2^{\binom{t}{2}} \le 2^{\binom{t}{2}+O(t \log t)}.$$
Using $w(S) \ge 1$ for every $S$, we finally have 
\begin{align*}
\|M(1_{\mathcal{B}} u)\|_w &\le \max_S |M(1_{\mathcal{B}} u)(S)|\\
&\le \|u\|_w \cdot 2^{\binom{t}{2}+O(t \log t)} \cdot 2^{-\binom{t}{2}-t\sqrt{t}+O(t\log t)} \cdot 2^{O(t \log t)}\\
&\le \|u\|_w \cdot 2^{-t\sqrt{t}+O(t \log t)}.
\end{align*}
In particular we can let $\gamma_3=2^{-ct^{3/2}}$ for some absolute constant $c>0$, such that $\|M(1_{\mathcal{B}} u)\|_w \le \gamma_3 \|u\|_w$ for large $t$.
\end{proof}

We wrap up all the previous estimates to give the missing proofs to earlier claims used in the proof of Theorem \ref{thm_small_n}.

\begin{theorem} \label{thm_final_estimate}
For sufficiently large $t$, and every $t \le n \le t^{16}$, for every $u \in \mathcal{C}$, we have 
$$\|T u\|_w=\|P 1_{\mathcal{B}} u\|_w \le \frac{1}{2}t^{-10} \|u\|_w.$$ 
With the $g_0$ we defined earlier, we also have 
$$\|P g_0\|_w \le \frac{1}{2}t^{-10}.$$
\end{theorem}
\begin{proof}
By Lemma \ref{lem_p_decompose}, we have $P=\sum_{W \subseteq [n]} A_W$. Lemma \ref{lem_small} tells us $A_W=0$ when $|W|<t$. Thus
$$P=\sum_{|W|=t} A_W + \sum_{t<|W|\le t+\ft} A_W + \sum_{|W|>t+\ft} A_W.$$
We use the estimates from Lemma \ref{lem_small}, \ref{lem_mid}, \ref{lem_large} for these three sums respectively. We get
$$\left\| P1_{\mathcal{B}} u\right\|_w \le (\gamma_1+\gamma_2+\gamma_3)\|u\|_w.$$
So it remains to check that $\gamma_1+\gamma_2+\gamma_3 \le \frac{1}{2}t^{-10}$ when $t$ is sufficiently large. It is easy to check that since $t \le n \le t^{16}$, we have
$$\gamma_1=4n^3 2^{-(t-1)/2} \le 4t^{48} 2^{-(t-1)/2} = o(t^{-10}),$$
$$\gamma_2=10tn^{5\sqrt{t}} 2^{-t} \sim 2^{c\sqrt{t}\log_2 t -t} =o(t^{-10}),\qquad \gamma_3=2^{-ct^{3/2}}=o(t^{-10}).$$
Adding the three bounds finishes the proof of $\|Tu\|_w \le \frac{1}{2} t^{-10} \|u\|_w$ for sufficiently large $t$.

For $g_0$, from Lemma \ref{lem_C} we know $g_0 \in \mathcal{C}$, this gives $$\|P 1_{\mathcal{B}} g_0\|_w=\|Tg_0\|_w \le \frac{1}{2}t^{-10}\|g_0\|_w.$$
In our definition of $g_0$, we set $g_0(S)=0$ whenever $S \not\in \mathcal{B}$. This means $1_{\mathcal{B}} g_0=g_0$. So we have $\|Pg_0\|_w \le \frac{1}{2}t^{-10}\|g_0\|_w$. By the definition of the $w$-norm and the choice of the weight function $w$,
$$\|g_0\|_w=\max_S \frac{|g_0(S)|}{w(S)}\le 1.$$
Therefore $$\|Pg_0\|_w \le \frac{1}{2}t^{-10},$$
completing the proof of the second assertion.
\end{proof}

\section{A technical interpolation lemma}\label{sec_interpolation}

We will extend Theorem \ref{thm_small_n} to larger $n>t^{16}$ in the next section. There we will take quotients of graphs on $n>t^{16}$ vertices with respect to random partitions of its vertex set into no more than $t^{16}$ parts, and apply $f_{t^{16}}$ given by Theorem \ref{thm_small_n} to define the satisfiable $f$ for Theorem \ref{thm_key} for large $n$. For that purpose, we rely on the following slightly technical lemma, which can be viewed as an interpolation result. Essentially one only needs to use the Lagrange interpolation theorem carefully, and estimate the sum of $|\beta_j|$. The proof is slightly tedious since $k_j$'s have to be integers.
\begin{lemma} \label{lem_interpolation}
Let $L \ge 3$ be an integer and $N=L^8$, set $$k_j=\left\lfloor \frac{N}{j^2} \right\rfloor \in [L^6, L^8], \qquad 1 \le j \le L.$$
Then there exist real numbers $\beta_1, \cdots, \beta_L$ such that
$$\sum_{j=1}^L \beta_j=1, \qquad \sum_{j=1}^L \beta_j k_j^{-d}=0\textup{~for~}1 \le d \le L-1, \qquad \sum_{j=1}^L |\beta_j| \le 4L.$$
\end{lemma}

\begin{proof}
First let us quickly check the rounding would not cause $k_j$ to be equal. For any $1 \le j \le L-1$,
$$\frac{N}{j^2}-\frac{N}{(j+1)^2}=N \cdot \frac{2j+1}{j^2(j+1)^2} \ge N \cdot \frac{2}{2L^4}=L^4 >3.$$
So the gap between two consecutive reals is at least $3$. After rounding, they are still pairwise distinct.

We apply the Lagrangian interpolation theorem. We set 
$$\beta_j=\prod_{1 \le i \le L: i \ne j} \frac{k_i^{-1}}{k_i^{-1}-k_j^{-1}}.$$
Then for every polynomial $p$ of degree at most $L-1$, we have
$$p(x)=\sum_{j=1}^{L} p(k_j^{-1}) \prod_{1 \le i \le L: i \neq j} \frac{x-k_i^{-1}}{k_j^{-1}-k_i^{-1}}.$$
Setting $x=0$ gives $$p(0)=\sum_{j=1}^L \beta_j p(k_j^{-1}).$$
Applying this identity to $p(x)=1, x, x^2, \cdots, x^{L-1}$ respectively, we get the first two assertions.

It remains to bound $\sum_{j=1}^L |\beta_j|$. We first ignore the rounding by taking $k'_j=\frac{N}{j^2}$ and let 
$$\beta'_j=\prod_{1 \le i \le L: i \ne j} \frac{(k'_i)^{-1}}{(k'_i)^{-1}-(k'_j)^{-1}}=\prod_{1 \le i \le L: i \neq j} \frac{i^2}{i^2-j^2}=2(-1)^{j-1} \cdot \frac{(L!)^2}{(L-j)! (L+j)!}.$$
Clearly $(L!)^2 \le (L-j)!(L+j)!$, so $|\beta_j'| \le 2$.

Now let us compare $\beta_j$ with $\beta'_j$. For each $i \neq j$, we have
$$\frac{k_i^{-1}}{(k'_i)^{-1}}=\frac{N}{i^2 \lfloor N/i^2 \rfloor} \le \frac{N}{N-i^2}\le \frac{N}{N-L^2} \le 1+\frac{2L^2}{N} \le \exp \left(\frac{2L^2}{N}\right).$$
The second last inequality follows from $L^2/N=1/L^6<1/2$. 

Note that 
$$(N/k_i)-i^2=\frac{N-i^2 k_i}{k_i}<\frac{i^2}{k_i}\le \frac{i^2}{N/(2i^2)}=\frac{2i^4}{N} \le \frac{2L^4}{N}.$$
Therefore we also have
\begin{align*}\left|\frac{(k'_i)^{-1}-(k'_j)^{-1}}{k_i^{-1}-k_j^{-1}}\right|&=\left|\frac{i^2-j^2}{N/k_i-N/k_j}\right|
\le \frac{1}{1-\frac{|(N/k_i-i^2)-(N/k_j-j^2)|}{|i^2-j^2|}}\\
& \le \frac{1}{1-2L^4/N} \le \exp\left(\frac{4L^4}{N}\right).
\end{align*}
Multiplying these two inequalities and then multiply for all $i \neq j$, we have
$$\left|\frac{\beta_j}{\beta'_j}\right|\le \exp\left(\frac{L(2L^2+4L^4)}{N}\right) \le \exp\left(\frac{6}{L^3}\right)<2.$$
The last inequality uses $L \ge 3$ and $e^{2/9}=1.248\cdots <2$. Therefore $|\beta_j|<2|\beta'_j| \le 4$. And thus 
$$\sum_{j=1}^L |\beta_j| \le 4L,$$
completing the proof of the last assertion.
\end{proof}

\section{Extending the function to every $n \ge t$}\label{sec_large_n}
The goal of this section is to prove Theorem \ref{thm_key} for sufficiently large $t$. For a positive integer $k$ and a graph $S$, consider a map $\phi: V(S) \rightarrow [k]$. The {\it quotient graph} $Q_\phi(S)$ is obtained from contracting vertices of $S$ that are mapped to the same label. In other words, $V(Q_\phi(S))$ has vertex set $[k]$, and two labels $i, j$ are adjacent if and only if there is one original edge $uv \in E(S)$ with $\phi(u)=i$, $\phi(v)=j$. 

Apply Lemma \ref{lem_interpolation} with $L=t^2$ (thus $N=t^{16})$, we have a sequence of real numbers $k_j$'s satisfying
$$\sum_{j=1}^L \beta_j=1, \qquad \sum_{j=1}^L \beta_j k_j^{-d}=0\textup{~for~}1 \le d \le L-1,$$
$$k_j = \left\lfloor \frac{t^{16}}{j^2} \right\rfloor \in [t^{12}, t^{16}],\qquad \sum_{j=1}^L |\beta_j| \le 4t^2.$$
We consider a random map $\phi: V(S) \rightarrow [k]$, such that every vertex is mapped to $[k]$ independently and uniformly at random. Define
\begin{equation}\label{eq_f}
f(S)=\sum_{j=1}^L \beta_j \cdot \mathbb{E}_{\phi: V(S) \rightarrow [k_j]} f_N(Q_{\phi}(S)).
\end{equation}
Here $f_N$ is the function $f$ given by Theorem \ref{thm_small_n} when $n=t^{16}=N$. After the contraction in the $j$-term, the quotient graph $Q_{\phi}(S)$ has $k_j \le t^{16}=N$ vertices. If needed, we may just append isolated vertices to it, so that $f_N$ can be evaluated properly. It is not hard to see that $f$ is isomorphism- and isolate-invariant. For the remainder of this section, we will verify that $f$ satisfies the rest of the requirements in Theorem \ref{thm_key}.

\begin{lemma}
The function $f$ defined in \eqref{eq_f} satisfies $f \in \mathcal{F}_{t-1}$. 
\end{lemma}
\begin{proof}
By the definition of $f$, it suffices to show that for every $n \ge t$, every $1 \le k \le N=t^{16}$, and every fixed map $\phi:[n] \rightarrow [k]$, the function $S \rightarrow f_N(Q_{\phi}(S))$ has Fourier support on $K_t$-free graphs.

To see this, first we observe that $f_N$ is Fourier supported on $K_t$-free graphs, meaning that 
$$f_N(T)=\sum_{H: H \subseteq K_N,\,\omega(H)<t} \widehat{f_N}(H)(-1)^{|H \cap T|}.$$
Substituting $T=Q_{\phi}(S)$, we have
$$f_N(Q_{\phi}(S))=\sum_{H: H \subseteq K_N,\,\omega(H)<t} \widehat{f_N}(H)(-1)^{|H \cap Q_{\phi}(S)|}.$$
We claim that $(-1)^{|H \cap Q_{\phi}(S)|}$ only depends on the following set $B$ of variables, represented by the edges they correspond to:
$$B=\{(u, v) \in E(K_n): \phi(u) \neq \phi(v), (\phi(u), \phi(v)) \in E(H)\}.$$
This can be derived from the simple observations: (i) for $\phi(u)=\phi(v)$, they do not give any edge in the quotient graph $Q_{\phi}(S)$;  (ii) for $(\phi(u), \phi(v)) \not\in E(H)$, whether $S$ contains that edge would not affect $H \cap Q_{\phi}(S)$.

For a function depending only a subset of variables, its Fourier coefficient is also supported in that set. Note that $\phi(B) \subseteq E(H)$ and the sum is over subgraphs $H$ of $K_N$ with $\omega(H)<t$. So the function $f_N(Q_{\phi}(S))$ is Fourier supported also on $K_t$-free graphs. This completes the proof.
\end{proof}

Next we show $f$ has the right value when the number of non-isolated vertices is at most $t$.
\begin{lemma}
For $f$ defined in \eqref{eq_f}, we have $f(\emptyset)=2^{\binom{t}{2}}-1$. Furthermore for $0<|V^+(S)| \le t$, we have $f(S)=(-1)^{e(S)+1}$, which gives $|f(S)| \le 1$.
\end{lemma}
\begin{proof}
Note that removing the isolated vertices from \eqref{eq_f} would not change the right hand side since $f_N$ was chosen to be isolate-invariant. So for fixed $k=k_j$, we have 
$$\mathbb{E}_{\phi: V(S) \rightarrow [k]} f_N(Q_{\phi}(S))=\mathbb{E}_{\phi: V^+(S) \rightarrow [k]} f_N(Q_{\phi}(S)).$$
Note that the quotient graph $Q_{\phi}(S)$ is determined by a partition $\pi$ of $V^+(S)$ into $k$ parts. We denote by $S/\pi$ the quotient graph with respect to $\pi$. Then
\begin{equation}\label{eq_quotient}
\mathbb{E}_{\phi: V^+(S) \rightarrow [k]} f_N(Q_{\phi}(S))=\sum_{\pi} f_N(S/\pi) \frac{k(k-1)\cdots (k-|\pi|+1)}{k^{|V^+(S)|}}.
\end{equation}
Here the sum is over all partition $\pi$ of $V^+(S)$ into $k$ parts, and $|\pi| \le k$ denotes the number of non-empty parts of $\pi$. 

Recall that $f_N$ was obtained from Theorem \ref{thm_small_n}. So $f_N(\emptyset)=2^{\binom{t}{2}}-1$ and $f_N(S)=(-1)^{e(S)+1}$ when $0<|V^+(S)| \le t$. For the quotient $S/\pi$, its number of non-isolated vertices is at most the number of non-empty parts $|\pi| \le |V^+(S)| \le t.$ So $f_N(S/\pi)$ equals $2^{\binom{t}{2}}-1$ if $S/\pi$ is empty, and equals $(-1)^{e(S/\pi)+1}$ otherwise.

Let us return to \eqref{eq_quotient}. When $S$ is the empty graph clearly $f_N(S/\pi)=2^{\binom{t}{2}}-1$ for every $\pi$ and thus the expectation is also equal to $2^{\binom{t}{2}}-1$. 

For the remainder of this proof, let us assume $0<|V^+(S)|\le t$. The right hand side of $\eqref{eq_quotient}$ can be viewed as a polynomial in $\frac{1}{k}$ of degree at most $|V^+(S)|-1$. We claim that its constant term is equal to $(-1)^{e(S)+1}$. Note that $k(k-1)\cdots (k-|\pi|+1)/k^{|V^+(S)|}$ contributes to the constant term if and only if $|\pi|=|V^+(S)|$, and it contributes exactly $1$ when that happens. Therefore the constant term on the right hand side is just equal to $f_N(S/\pi)=f_N(S)=(-1)^{e(S)+1}.$ Now let us assume that for every $k$,
$$\mathbb{E}_{\phi: V^+(S) \rightarrow [k]} f_N(Q_{\phi}(S))=(-1)^{e(S)+1}+\sum_{r=1}^{|V^+(S)|-1}\frac{c_{r}}{k^r}.$$
Then by definition,
\begin{align*}
f(S)&=\sum_{j=1}^L \beta_j \cdot \mathbb{E}_{\phi: V^+(S) \rightarrow [k_j]} f_N(Q_{\phi}(S))\\
&=\sum_{j=1}^L \beta_j\left((-1)^{e(S)+1}+\sum_{r=1}^{|V^+(S)|-1}\frac{c_{r}}{k_j^r}\right)\\
&=(-1)^{e(S)+1} \sum_{j=1}^L \beta_j + \sum_{r=1}^{|V^+(S)|-1}   c_{r}\left(\sum_{j=1}^L \beta_j k_j^{-r}\right)\\
&=(-1)^{e(S)+1}.
\end{align*}
The last inequality follows from $\sum_{j=1}^L \beta_j=0$ and $\sum_{j=1}^L \beta_j k_j^{-r}=0$ for $1\le r \le |V^+(S)|-1 \le t-1 \le L-1.$
\end{proof}
Finally, we verify (3) in Theorem \ref{thm_key}. 
\begin{lemma}
For $f$ defined in \eqref{eq_f} and $|V^+(S)|>t$, we have
$$|f(S)| \le t^{-7},$$
when $t$ is sufficiently large.
\end{lemma}
\begin{proof}
In the sum \eqref{eq_f}, we consider three types of quotient graphs $Q_{\phi}(S)$: empty quotients, non-empty quotients with at most $t$ non-isolated vertices, and quotients with more than $t$ non-isolated vertices. 

We start with the second class with $0<|V^+(Q_{\phi}(S))| \le t$. We choose edges of $S$ in an arbitrary order, add them one at a time until we reach its first subgraph $S_0$ with more than $t$ non-isolated vertices. This step does not involve $\phi$. Such subgraph exists since we start with zero non-isolated vertices, and end with $S$ that has more than $t$ non-isolated vertices. Since adding each edge increase $|V^+(\cdot)|$ by at most $2$. We know $t+1 \le |V^+(S_0)| \le t+2$. Suppose $\phi$ is injective on $V^+(S_0)$, then 
$$|V^+(Q_{\phi}(S))| \ge |V^+(Q_{\phi}(S_0))| \ge |V^+(S_0)| \ge t+1,$$
contradicting $|V^+(Q_{\phi}(S))| \le t$. So we know that $\phi$ cannot be injective on $V^+(S_0)$. For fixed $k_j$, the probability that $\phi$ is not injective is at most $\binom{|V^+(S_0)|}{2} \cdot \frac{1}{k_j} \le \binom{t+2}{2} \cdot \frac{1}{k_j}$. This quantity also upper bounds the probability that a random $\phi$ gives a quotient belonging to the second class. Note that for such quotient, from the proof of previous lemma we know $|f_N(Q_{\phi}(S))|=1$. 

For the first class of quotient graph $Q_{\phi}(S)=\emptyset$, obviously $\phi$ is constant on each connected component of $S$. Let $c(S)$ be the number of non-empty connected components of $S$, then
$$\mathbb{P}_{\phi: V(S) \rightarrow [k_j]}[Q_{\phi}(S)=\emptyset]=k_j^{c(S)-|V^+(S)|}.$$
For such $\phi$, $f_N(Q_{\phi}(S))=2^{\binom{t}{2}}-1.$

Finally, for the third class of quotient graphs when $|V^+(Q_{\phi}(S))|>t$, Theorem \ref{thm_small_n} gives $|f_N(Q_{\phi}(S))| \le t^{-10}$. Combining all these three estimates, we have
\begin{align*}
f(S)&=\sum_{j=1}^L \beta_j \cdot \mathbb{E}_{\phi: V(S) \rightarrow [k_j]} f_N(Q_{\phi}(S))\\
& \le (2^{\binom{t}{2}}-1)\sum_{j=1}^L \beta_j k_j^{-(|V^+(S)|-c(S))} + \sum_{j=1}^{L} \binom{t+2}{2}\cdot \frac{|\beta_j|}{k_j}+ t^{-10} \cdot \sum_{j=1}^{L} |\beta_j|.
\end{align*}
For the first sum, we have $|V^+(S)|-c(S) \ge 1$ since $S$ is non-empty. Using $\sum_{j=1}^L \beta_j k_j^{-d}=0$ for $1 \le d \le L-1$, this sum is either $0$, or at most $(2^{\binom{t}{2}}-1)t^{-12L} \sum_{j=1}^L |\beta_j|=(2^{\binom{t}{2}}-1)t^{-12t^2} \sum_{j=1}^L |\beta_j|$, since only those terms with $|V^+(S)-c(S)| \ge L$ will remain and $k_j \ge t^{12}$ from Lemma \ref{lem_interpolation}. Therefore using $\sum_{j=1}^{L} |\beta_j| \le 4t^2$, we have
\begin{align*}
|f(S)| &\le 2^{t^2/2} \cdot 4t^2 \cdot t^{-12t^2} + t^2 \cdot 4t^2 \cdot t^{-12}+ t^{-10}\cdot 4t^2=O(t^{-8}) \le t^{-7}
\end{align*}
when $t$ is sufficiently large.

\end{proof}

\section{A new proof for the odd-cycle-intersecting case}\label{sec_k3}
A family $\mathcal{F}$ of graphs is called {\it odd-cycle-intersecting}, if the intersection of two graphs $G_1, G_2 \in  \mathcal{F}$ contains an odd cycle, i.e. $\chi(G_1 \cap G_2) \ge 3$. Clearly a $K_3$-intersecting family is also odd-cycle-intersecting, but the converse might not be true.

In this section, we give a new proof for the triangle-intersecting case by showing that every odd-cycle-intersecting family has density at most $\frac{1}{8}$. Our proof is slightly different from the proof by Ellis, Filmus, and Friedgut in \cite{EFF12} and the one by Berger and Zhao \cite{BZ23}. The main new ingredient in this proof is the monotonicity of the maximum family. 

\begin{theorem}
For $n \ge 3$, and an odd-cycle-intersecting family $\mathcal{F}$ of graphs on vertex set $[n]$, we always have 
$$|\mathcal{F}| \le 2^{\binom{n}{2}-3}.$$
Moreover, when equality is attained, $\mathcal{F}$ consists of graphs on $[n]$ that contains a fixed copy of triangle.
\end{theorem}

\begin{proof}
Take an odd-cycle-intersecting family $\mathcal{F}$ of graphs on $[n]$. We view each graph in $\mathcal{F}$ as a $\{-1,1\}$-vector $x=(x_e)_{e \in E(K_n)}$ in $\mathbb{R}^{E(K_n)}$, with $+1$ denoting a present edge. Consider the indicator function $f=1_{\mathcal{F}}: \{-1,1\}^{E(K_n)} \rightarrow \mathbb{R}$ which outputs $1$ if the input vector corresponds to a graph from $\mathcal{F}$, and outputs $0$ otherwise. By Fourier transform with respect to the standard basis $\prod_{e \in S} x_e$ for all $S \subseteq E(K_n)$, we have:
$$f(x)=\sum_{S \subseteq E(K_n)} \widehat{f}(S) \chi_S(x),$$
where $\chi_S(x)=\prod_{e \in S} x_e.$

Since $\mathcal{F}$ is odd-cycle-intersecting, in particular it is intersecting and thus $f(x)f(-x)=0$. Therefore if we let $h(x)=\frac{f(x)-f(-x)}{2}$ and $\mu=|\mathcal{F}|/2^{\binom{n}{2}}$, then we have $\widehat{f}(\emptyset)=\mu$, and furthermore
$$\|h\|_2^2=\mathbb{E}\left[\frac{(f(x)-f(-x))^2}{4}\right]=\frac{\mathbb{E}[f(x)^2+f(-x)^2]}{4}=\frac{\mu}{2}.$$
Also $\widehat{h}(S)=\widehat{f}(S)$ when $|S|$ is odd and $0$ when $|S|$ is even.

First let us observe that $f(x)f(\overline{x \oplus C})=0$ for any cut $C=E(V_1, V_2)$. In other words, if we take the complement of the graph obtained from switching the adjacency between vertices in $V_1$ and $V_2$ in $x$, then it cannot co-exist with $x$ in $\mathcal{F}$, this is because as graphs 
$$x \cap (\overline{x \oplus C})=x \cap C.$$
Since $C$ is bipartite, the intersection is also bipartite and thus cannot contain an odd cycle. Now we take a random cut $C$ by assigning every vertex of $K_n$ independently with probability $1/2$ each into $V_1$ and $V_2$. For each such $C$, we can view it as assigning $\sigma_v \in \{-1, 1\}$ to each vertex of $K_n$, and for the edge $e=uv$, $x_{uv}$ becomes $-x_{uv}\sigma_u \sigma_v$ in $\overline{x \oplus C}.$
As a result,
\begin{align*}
0&=\mathbb{E}_{C,x}[f(x)f(\overline{x \oplus C})]\\
&=\mathbb{E}_{C,x} \left[\left(\sum_S \widehat{f}(S) \chi_S(x)\right)\left(\sum_S \widehat{f}(S) \chi_S(\overline{x \oplus C})\right)\right]\\
&=\mathbb{E}_\sigma\left[\sum_S \widehat{f}(S)^2 \cdot (-1)^{|S|} \cdot \prod_{uv \in S} \sigma_u \sigma_v\right]
\end{align*}
For a subgraph $S$ of $K_n$, let $O(S)$ be the set of odd-degree vertices in $S$. Note that $\mathbb{E}_{\sigma_u} \sigma_u^d=0$ for odd $d$ and $1$ for even $d$. Therefore $\mathbb{E}_\sigma [\prod_{uv \in S} \sigma_u \sigma_v]$ equals $1$ if and only if $O(S)=\emptyset$, and equals $0$ otherwise. Therefore, 
$$0=\mathbb{E}_C[f(x)f(\overline{x \oplus C})]=\sum_{S: O(S)=\emptyset} (-1)^{|S|} \widehat{f}(S)^2.$$
This gives
\begin{equation}\label{eqk3_1}
\sum_{S:O(S)=\emptyset} \widehat{h}(S)^2=\sum_{S:O(S)=\emptyset,|S|\textup{~odd}} \widehat{f}(S)^2=\sum_{S:O(S)=\emptyset,|S|\textup{~even}} \widehat{f}(S)^2 \ge \widehat{f}(\emptyset)^2 = \mu^2.
\end{equation}

We may assume that $f$ is monotone increasing, since replacing $\mathcal{F}$ by its upward closure preserves the hypothesis and can only increase its size. So $h(x)=(f(x)-f(-x))/2$ is also monotone increasing. Now we take the derivative along the edge $e$. Let $y$ be the vector obtained from removing the $e$-th coordinate of $x$, and define
$$D_e h(y):=\frac{h(x_e\leftarrow 1, y)-h(x_e\leftarrow -1, y)}{2}.$$
Then clearly $D_e h(y) \ge 0$ from the monotonicity of $h$. Also its Fourier expansion can be written as 
$$D_e h(y)=\sum_{S: e \in S} \widehat{h}(S) \chi_{S\setminus e}(y).$$
We take a random cut $C$ such that the edge $e \in E(C)$. If we assume $e=uv$, then this corresponds to a random assignment of $\{-1, 1\}$ to $\{\sigma_v\}_{v \in V(K_n)}$ such that $\sigma_u\sigma_v=-1$. We have $D_e h(y)D_e h(y \oplus C) \ge 0$. Taking the expectation, this gives
\begin{align*}
0 &\le \mathbb{E}_C [D_e h(y) D_e h(y \oplus C)]\\
&=\mathbb{E}_C\left[\left(\sum_{S: e \in S}\widehat{h}(S) \chi_{S\setminus e}(y)\right)\left(\sum_{S: e \in S}\widehat{h}(S) \chi_{S\setminus e}(y \oplus C)\right)\right]\\
&=\mathbb{E}_\sigma \left[\sum_{S: e \in S}\widehat{h}(S)^2 \prod_{uv \in S \setminus e} \sigma_u \sigma_v \right]
\end{align*}
Since $\sigma$ is uniformly random over those with $\sigma_u \sigma_v=-1$, it is not hard to see that $\mathbb{E}_{\sigma}[\prod_{uv \in S\setminus e} \sigma_u \sigma_v]$ equals $-1$ if $O(S)=\emptyset$, $1$ if $O(S)=\{u, v\}$, and $0$ in all other cases. Therefore, the inequality above becomes:
$$\sum_{S: e \in S, O(S)=\emptyset} \widehat{h}(S)^2 \le \sum_{S: e \in S, O(S)=\{u,v\}=e} \widehat{h}(S)^2.$$
Summing over all edges $e \in E(K_n)$, we have
\begin{equation}\label{eq_k3_2}
\sum_{S: O(S)=\emptyset} |S| \widehat{h}(S)^2 \le \sum_{S: O(S) \neq \emptyset} \widehat{h}(S)^2.
\end{equation}
This comes from the fact that every $S$ with $O(S)=\{u, v\}$ appears at most once when summing the right hand side. 

Recall that $\widehat{h}(S)$=0 when $|S|$ is even, and when $O(S)=\emptyset$ and $S \neq \emptyset$, $S$ has at least $3$ edges. Therefore the left hand side is at least $3 \sum_{S: O(S)=\emptyset} \widehat{h}(S)^2$. So we have
$$3 \sum_{S: O(S)=\emptyset} \widehat{h}(S)^2 \le \sum_{S: O(S)\neq \emptyset} \widehat{h}(S)^2.$$
Using the fact that
$$\sum_{S: O(S)=\emptyset} \widehat{h}(S)^2 + \sum_{S: O(S)\neq\emptyset} \widehat{h}(S)^2=\mathbb{E}[h^2]=\frac{\mu}{2}.$$
We immediately have 
$$\sum_{S: O(S)=\emptyset} \widehat{h}(S)^2 \le \frac{\mu}{8}.$$
This gives $\mu^2 \le \mu/8$, and $\mu \le \frac{1}{8}$, and the proof of the upper bound is complete. 

Next we prove the uniqueness part. Observe that when the equality is attained, from \eqref{eqk3_1}, $\widehat{f}(S)=0$ whenever $|S|>0$ is even and $O(S)=\emptyset$. From how we bound $|S|$ by $3$ in \eqref{eq_k3_2}, we also know that $\widehat{f}(S)=0$, whenever $|S| > 3$ is odd and $O(S)=\emptyset$. Therefore the Fourier coefficient of $f$ when $S$ is Eulerian is supported on the empty graph and triangles. Using the same averaging argument before, we have that for the random cut $C$ (which is bipartite and thus not in $\mathcal{F}$),
$$0=\mathbb{E}_C[f(C)]=\mathbb{E}_C \left[\sum_S \widehat{f}(S)\chi_S(C) \right]=\sum_{S: O(S)=\emptyset} (-1)^{|S|}\widehat{f}(S)=\mu-\sum_{T~\textup{triangle}}\widehat{f}(T),$$
where the last sum is over all triangles $T$. This immediately gives $\sum_{T\textup{~triangle}} \widehat{f}(T)=\mu$.

For each edge $e$, we let $p_e=\mathbb{P}_C[C \oplus \{e\} \in \mathcal{F}]$. Note that for a triangle $T$, a cut $C$ meets $T$ in either two or zero edges, therefore $\chi_T(C)$ is always equal to $-1$. We can also see that $\chi_T(C \oplus \{e\})=1$ if $e \in T$ and $-1$ if $e \not\in T$. Therefore,
\begin{align*}
p_e&=\mathbb{E}_C f(C\oplus \{e\})=\mathbb{E}_C \left[\sum_S \widehat{f}(S) \chi_S(C \oplus \{e\})\right]\\
&=\mu+\left(\sum_{T~\textup{triangle}: e \in T} \widehat{f}(T)\right) - \left(\sum_{T~\textup{triangle}: e \not\in T} \widehat{f}(T)\right)=2 \sum_{T~\textup{triangle}: e \in T} \widehat{f}(T).
\end{align*}
Summing this over all $e$ gives
$$\sum_e p_e = 6 \sum_{T\textup{~triangle}} \widehat{f}(T)=6\mu.$$
One can similarly get $0 \le \mathbb{E}_C f(\overline{C \oplus \{e\}}) = 2 \mu-p_e$. So we know all $p_e$ lie in the interval $[0, 2\mu]$. Since their sum is $6\mu$, there are at least three edges $e_1, e_2, e_3$ such that for each $i=1,2,3$, $p_{e_i}>0$, meaning that there exists some cut $C$ such that $C \oplus \{e_i\} \in \mathcal{F}$, clearly $C$ cannot contain $e_i$, otherwise $C \oplus \{e_i\}=C \setminus \{e_i\}$ which becomes bipartite. On the other hand, for any other graph $G \in \mathcal{F}$, if $e_i \not\in E(G)$, $G \cap (C \oplus e_i)$ is a subgraph of $C$ which is again bipartite. So every graph in $\mathcal{F}$ must contain the edge $e_i$ for each of $i=1,2,3$. This already gives $\mathcal{F}$ of density at most $\frac{1}{8}$. When equality holds, to make it odd-cycle-intersecting, we clearly need the three edges $e_1,e_2,e_3$ to form a triangle. This completes the proof of the uniqueness.
\end{proof}

\section{Concluding remarks}\label{sec_concluding}
In this paper we show a $K_t$-intersecting family of graphs has density at most $2^{-\binom{t}{2}}$ for sufficiently large $t$. The main obstruction for our approach to work for arbitrary $t$ is the way of constructing functions $g$ (that gives $f$) in the proof of Theorem \ref{thm_small_n} stops working once $n$ is large in $t$. One possibility is to replace $T=P1_{\mathcal{B}}$ by $PQ$ for some other function $Q$ that is supported on $\mathcal{B}$. Also the formula we used for $g$ is the same for all $t \le n \le t^{16}$. To make the construction work for all $t$, one would probably need to define $g$ inductively as the number of vertices $n$ grows.

Filmus \cite{F11} gave an elementary proof for the triangle-intersecting result for $n=8$. A more combinatorial way to understand why having eight vertices is special is to take $H=\overline{Q^3}$, where $Q^3$ is the $3$-dimensional hypercube graph, then for every $8$-vertex graph $G$, $\mathbb{E}_{\pi \in S_8} (-1)^{e(G \cap \pi(H))} \ge -\frac{1}{7}$. This lower bound can be obtained from the fact that as a $3$-uniform hypergraph, the Fano plane has $7$ vertices and independence number equal to $4$. One might hope that for arbitrary $t$, one can find constantly many (say one or two) $n$-vertex graphs $H_1, \cdots, H_c$ whose complements are $K_t$-free, and some linear combination of $\mathbb{E}_{\pi \in S_n} (-1)^{e(G \cap \pi(H))}$ that has maximum $1$ and is always above $-1/(2^{\binom{t}{2}}-1)$. Suppose this works for $n$ exponential in $t$, then perhaps the interpolation trick in Section \ref{sec_interpolation} and \ref{sec_large_n} can be used to extend the function $f$ to all $n$.

For general $H$, the following question remains open.

\begin{question}
For a connected graph $H$, is it true that the maximum density of an $H$-intersecting family of graphs on vertex set $[n]$ for every $n \ge |V(H)|$ is precisely $2^{-e(H)}$, if and only if $H$ is vertex-transitive?
\end{question}
This characterization is likely far from correct, yet, somewhat surprisingly,  we know of no counterexample in either direction. It would be interesting to determine whether $H=C_4$ has this property, or to construct a connected and non-vertex-transitive graph that does.

The following question appears more tractable compared to the problems discussed above, and we expect an affirmative
answer for all $n \ge t > k$. For $k \ge 3$, the only case
settled so far is that of the tetrahedron $K_4^3$, established
by Mani and Zhang \cite{MZ25}.
\begin{question}
Can we apply the same technique to prove that for fixed $k$ and sufficiently large $t > k$ and every $n \ge t$, every $K_t^k$-intersecting family of $k$-uniform hypergraphs on vertex set $[n]$ has size at most $2^{\binom{n}{k}-\binom{t}{k}}$?
\end{question}

One might also wonder, can the method introduced in this paper be applied to other kinds of intersection problems? Here is one problem for which we believe our technique might be useful. Take a subset $B$ of the cyclic group $\mathbb{Z}/(n\mathbb{Z})$, a family $\mathcal{F}$ of subsets of $\mathbb{Z}/(n\mathbb{Z})$ is called {\it $B$-translate-intersecting} if for every $S_1, S_2 \in \mathcal{F}$, $S_1 \cap S_2$ contains a translate of $B$ in the group. Chung, Frankl, Graham and Shearer \cite{CFGS86} conjectured that the maximum density of $\mathcal{F}$ is $2^{-|B|}$ and they proved this when $B$ is an interval. Russell \cite{R09} gave a different proof. F\"uredi, Griggs, Holzman, and Kleitman \cite{FGHK90} proved this conjecture for all $|B|=3$. Griggs and Walker \cite{GW89} proved it for fixed $B$ and infinitely many $n$. Perhaps our result, in particular the interpolation part, can be used for sufficiently large $|B|$ and arbitrary $n \ge |B|$, at least for certain well-structured classes beyond the intervals.

We conclude the paper with one more question, inspired by the recent proofs and subsequent strengthenings of the famous conjecture of Chv\'atal \cite{Chv74} by Chang, Liu, and Liu \cite{CLL26}, Keevash \cite{K26}, and Ellis, Filmus, and Friedgut \cite{EFF26}. A family $\mathcal{G}$ of graphs on vertex set $[n]$ is {\it hereditary}, if for any $G_1 \subseteq G_2$, $G_2 \in \mathcal{G}$, we have $G_1 \in \mathcal{G}$.
\begin{question}
Is it true that for every $t \ge 3$, and every hereditary family $\mathcal{G}$ of graphs on $[n]$ with $n \ge t$, the maximum size of its $K_t$-intersecting subfamily is attained by taking all graphs from $\mathcal{G}$ containing some fixed copy of $K_t$?
\end{question}
This question seeks a common generalization of the $K_t$-intersecting conjecture and now-resolved Chv\'atal's Conjecture. For $t=2$ this follows directly from Chv\'atal. It would be interesting to see if techniques developed in these recent
works can be combined with our argument in
Section \ref{sec_k3} to give a positive answer to this question for $t=3$.\\~\\

\noindent {\bf Declaration of AI use.~} Multiple versions of ChatGPT (mainly 5.6 Sol and 6 Astra) were used by the authors over the past several months to test their ideas, write codes to solve numerical linear programs, and produce the technical estimates. All AI-assisted outputs were independently checked, verified, and simplified by the authors, who assume full responsibility for the content of this work.\\
\appendix

\section{Proof of the technical Lemma \ref{lem_pair}}
\begin{proof}
First we prove that when $t$ is sufficiently large and $n\in[t,t^{16}]$, for every integer $r \ge 0$ and every graph $S$,
$$\frac{1}{w(S)} \sum_{W \subseteq [n]: |W|=t+r} \sum_{T \in \mathcal{B}: S \Delta T \subseteq K_W} w(T) \le 4n^{2r}2^{\binom{t}{2}}.$$
We can assume $n\ge t+r$ since otherwise the inequality is trivially true. Write $q=n^3 2^{-(t-1)/2}$. We may assume $q<1/3$ since $q\le t^{48}2^{-(t-1)/2}\to 0$. We further observe that $q<1/3$  also implies $2^{\binom{t}{2}}\ge n^t$. Recall that 
$$w(S)=\begin{cases}
2^{\binom{t}{2}}, \qquad & |V^+(S)|=0,\\
n^{t-|V^+(S)|},\qquad & 1\le|V^+(S)|\le t,\\
1,\qquad & |V^+(S)|>t.
\end{cases}$$
Notice that in all cases $w(S)\ge n^{t-|V^+(S)|}$. Write $k=|V^+(S)|$, we shall replace $w(S)$ (but not $w(T)$ inside the sum) by $n^{t-k}$ and show that the inequality is still true.
Since $w(T)$ depends only on $|V^+(T)|$, our strategy is to group the sum by $|V^+(T)|$ and turn it into a counting problem. Let
$$
    c_j = \#\{W\subseteq[n], T\subseteq K_n: |W|=t+r, |V^+(T)|=j, S\Delta T\subseteq K_W\}
$$
for $0\le j\le t$. Then the left hand side becomes
$$
\frac{1}{n^{t-k}}\left(2^{\binom{t}{2}}c_0+\sum_{j=1}^t n^{t-j}c_j\right).
$$
When $j=0$, the only graph $T$ with $|V^+(T)|=0$ is $T=\emptyset$, thus the condition $S\Delta T\subseteq K_W$ becomes $S\subseteq K_W$, so $c_0$ is simply
$$\#\{W \subseteq [n]: |W|=t+r, S  \subseteq K_W\}=\binom{n-k}{t+r-k}\le n^{t+r-k}.$$
For $j\ge 1$, the condition $S\Delta T\subseteq K_W$ means that $S$ and $T$ have the same non-isolated vertices outside $W$. So we can estimate the number of $W$ and $T$ for fixed $|V^+(S)\setminus W|=\ell$ and sum them all together. For $W$, we shall choose $k-\ell$ vertices from $V^+(S)$ and $t+r+\ell-k$ vertices from $[n]\setminus V^+(S)$, which gives us
$$
\binom{k}{\ell}\binom{n-k}{t+r+\ell-k}\le n^{t+r+2\ell-k}
$$
choices.
For $T$, we shall consider all subgraphs of $K_W$ with $|V^+(T)|-|V^+(T)\setminus W|=j-\ell$ vertices, ignoring whether it makes a valid $T$, which may only increase the upper bound. This gives us at most 
$$
\binom{t+r}{j-\ell}2^{\binom{j-\ell}{2}}=\binom{t+r}{t+r-j+\ell}2^{\binom{j-\ell}{2}}\le n^{t+r-j+\ell}2^{\binom{j-\ell}{2}}
$$
choices. Summing over all $\ell$ from $0$ to $t$, we get

$$
    c_j\le\sum_{\ell=0}^{j} n^{t+r+2\ell-k}n^{t+r-j+\ell}2^{\binom{j-\ell}{2}}=n^{2r}\sum_{\ell=0}^{j} n^{2t+3\ell-k-j}2^{\binom{j-\ell}{2}}.
$$
Plug the upper bounds for $c_j$'s in and clean it up, it's enough to show that 
$$
2^{\binom{t}{2}}n^r+n^{2r}\sum_{j=1}^t \sum_{\ell=0}^{j} n^{2t+3\ell-2j}2^{\binom{j-\ell}{2}}\le 4n^{2r}2^{\binom{t}{2}}.
$$
Write $m=t+\ell-j$, then $j-\ell=t-m, \ell\le m$ and we have
$$
    \sum_{j=1}^t \sum_{\ell=0}^{j} n^{2t+3\ell-2j}2^{\binom{j-\ell}{2}}\le \sum_{m=0}^{t}(m+1)n^{3m}2^{\binom{t-m}{2}}.
$$
Using $\binom{t}{2}-\binom{t-m}{2}\ge t(t-1)/2-(t-m)(t-1)/2=m(t-1)/2$ and $q<1/3$, we have 
$$
 \sum_{m=0}^{t}(m+1)n^{2m}2^{\binom{t-m}{2}}\le  2^{\binom{t}{2}}\sum_{m=0}^{t}(m+1)n^{2m}2^{-m(t-1)/2}=2^{\binom{t}{2}}\sum_{m=0}^{t}(m+1)q^m\le\frac{2^{\binom{t}{2}}}{(1-q)^2}\le 3 \cdot 2^{\binom{t}{2}}
$$
and we are done. 

Then we are going to prove the second part, that is,
$$\frac{1}{w(S)} \sum_{W \subseteq [n]: |W|=t, S\nsubseteq K_W} \sum_{T \in \mathcal{B}: S \Delta T \subseteq K_W} w(T) \le 4n^{3}2^{(t-1)^2/2}.$$
Similarly set
$$
    c'_j= \#\{W\subseteq[n], T\subseteq K_n: |W|=t, |V^+(T)|=j, S\Delta T\subseteq K_W, S\nsubseteq K_W\}.
$$
Then 
$$
c'_0=\#\{W\subseteq[n]: |W|=t, S\subseteq K_W, S\nsubseteq K_W\}=0.
$$
When $j\ge1$, the estimate is similar except that $\ell=|V^+(S)\setminus W|\ge 1$ and $r=0$. Still write $m=t+\ell-j$, the upper bound becomes
$$
2^{\binom{t}{2}}\sum_{m=1}^{t}m q^m\le2^{\binom{t}{2}}\frac{q}{(1-q)^2}\le 4 \cdot 2^{\binom{t}{2}}q=4n^32^{(t-1)^2/2},
$$
completing the proof of the second part of this lemma.
\end{proof}

\end{document}